\documentclass[3p, final]{elsarticle}
\usepackage{amssymb}
 \usepackage{amsthm}
\usepackage{amsmath,bbm,dsfont, color, xfrac, accents}

\usepackage[OT2, T1]{fontenc}

\usepackage{showkeys}

\journal{arXiv}

\begin{document}
\newcommand{\csubset}         {\hookrightarrow}
\newcommand*{\id}{\mathds{1}}
\newcommand*{\zahlen}{\mathbb{Z}}
\newcommand*{\en}{\mathbb{N}}
\newcommand*{\er}{\mathbb{R}}
\newcommand{\divv}{\mathop{\text{div}}}
\newcommand{\esssup}{\mathop{\text{ess\:sup}\,}}
\newcommand{\htimes}{\mathop{\text{\large$�$}}}
\newcommand{\hexists}{\mathop{\text{\LARGE$\exists$}}}
\newcommand{\hforall}{\mathop{\text{\LARGE$\forall$}}}
\newtheorem{lem}{Lemma}
\newtheorem{teo}{Theorem}
\newtheorem{cor}{Corollary}
\newtheorem{prop}{Proposition}
\newtheorem{mydef}{Definition}
\newtheorem{rem}{Remark}
\newtheorem{ass}{Assumption}

\newcommand{\eps} {\varepsilon}


\newcommand{\Div}       {{\rm div}_x}
\newcommand{\toc}       {{\stackrel{b}{\longrightarrow\,}}}

\newcommand{\A}     {{\mathcal A}}

\newcommand{\Cam}     {{\mathcal L}}

\newcommand{\V}     {{\mathcal V}}

\newcommand{\D}     {{\tens {D}}}
\newcommand{\Dd}     {{\tens {D^d}}}
\newcommand{\nablat}     {{\tens {\nabla}}}

\def\tens#1{\pmb{\mathsf{#1}}}
\def\vec#1{\boldsymbol{#1}}

\newcommand{\eN}	{{\mathcal{N}}}
\newcommand{\M}	{{\mathcal{M}}}
\newcommand{\eO}	{{\mathcal{O}}}
\newcommand{\T}	{{\mathcal{T}}}
\newcommand{\W}	{\mathbb W}
\newcommand{\ti}	{\tilde}
\newcommand{\wti}	{\widetilde}
\newcommand{\vep}{\varepsilon}
\newcommand{\Divv}{{\rm {\bf div \,}}}
\newcommand{\Om}{\Omega}

\newcommand{\vdd}{\tens{\bar{v}}}
\newcommand{\Ref} {\eqref}

\newcommand{\jb}[1]{\textcolor{red}{\sc[**JB: #1 **]}}

\newcommand{\ddd}{\,{\rm d}}

\def\note#1{\marginpar{\small #1}}

\def\tens#1{\pmb{\mathsf{#1}}}

\def\vec#1{\boldsymbol{#1}}

\def\norm#1{\left|\!\left| #1 \right|\!\right|}

\def\fnorm#1{|\!| #1 |\!|}

\def\abs#1{\left| #1 \right|}

\def\ti{\text{I}}

\def\tii{\text{I\!I}}

\def\tiii{\text{I\!I\!I}}

\newcommand{\supp}{\operatorname{supp}}
\newcommand{\esup}{\operatorname{ess\,sup\,}\displaylimits}

\def\diver{\mathop{\mathrm{div}}\nolimits}

\def\grad{\mathop{\mathrm{grad}}\nolimits}

\def\Div{\mathop{\mathrm{Div}}\nolimits}

\def\Grad{\mathop{\mathrm{Grad}}\nolimits}

\def\tr{\mathop{\mathrm{tr}}\nolimits}

\def\cof{\mathop{\mathrm{cof}}\nolimits}

\def\det{\mathop{\mathrm{det}}\nolimits}

\def\lin{\mathop{\mathrm{span}}\nolimits}

\def\pr{\noindent \textbf{Proof: }}

\def\pp#1#2{\frac{\partial #1}{\partial #2}}

\def\dd#1#2{\frac{\d #1}{\d #2}}

\def\bA{\tens{A}}

\def\T{\mathcal{T}}

\def\R{\mathcal{R}}

\def\bx{\vec{x}}

\def\be{\vec{e}}

\def\bef{{f}}

\def\bec{\vec{c}}

\def\bs{\vec{s}}

\def\ba{\vec{a}}

\def\bn{{n}}

\def\bphi{\vec{\varphi}}

\def\btau{\vec{\tau}}

\def\bc{\vec{c}}

\def\bg{{g}}

\def\mO{\mathcal{O}}
\def\pmO{\partial\mathcal{O}}

\def\bE{\tens{\varepsilon}}
\def\bsig{\tens{\sigma}}

\def\bW{\tens{W}}

\def\bT{\tens{T}}

\def\bxi{\tens{\xi}}

\def\bD{\tens{D}}

\def\bF{\tens{F}}

\def\bB{\tens{B}}

\def\bV{\tens{V}}

\def\bS{\tens{S}}

\def\bI{\tens{I}}

\def\bi{\vec{i}}

\def\bv{{v}}

\def\bfi{\vec{\varphi}}

\def\bk{\vec{k}}

\def\b0{\vec{0}}

\def\bom{\vec{\omega}}

\def\bw{{w}}

\def\p{\pi}

\def\bu{{u}}
\def\bz{\vec{z}}
\def\bep{\vec{e}_{\textrm{p}}}
\def\dbep{\dot{\vec{e}}_{\textrm{p}}}
\def\bee{\vec{e}_{\textrm{el}}}
\def\dbee{\dot{\vec{e}}_{\textrm{el}}}

\def\ID{\mathcal{I}_{\bD}}

\def\IP{\mathcal{I}_{p}}

\def\Pn{(\mathcal{P})}

\def\Pe{(\mathcal{P}^{\eta})}

\def\Pee{(\mathcal{P}^{\varepsilon, \eta})}

\def\dx{\,{\rm dx}}

\begin{frontmatter}

 \title{Existence and smoothness for a class of $n$D models in elasticity theory of small deformations}

\author{Miroslav Bul\'{\i}\v{c}ek}  \address{Mathematical Institute, Charles University,  Sokolovsk\'a 83, 186 75 Prague 8, Czech Republic}

\author{Jan Burczak} \address{Institute of Mathematics, Polish Academy of Sciences, \'Sniadeckich 8, 00-656 Warsaw, Poland \\
OxPDE, Mathematical Institute, University of Oxford, Oxford, UK}



\begin{abstract}
We consider a model for deformations of a homogeneous isotropic body, whose shear modulus remains constant, but its bulk modulus can be a highly nonlinear function. We show that for a general class of such models, in an arbitrary space dimension, the respective PDE problem has a unique solution. Moreover, this solution enjoys interior smoothness. This is the first regularity result for elasticity problems that covers the most natural space dimension $3$ and that captures behaviour of many typical elastic materials (considered in the small deformations) like rubber, polymer gels or concrete.
\end{abstract}




\end{frontmatter}





\newenvironment{rcases}
  {\left.\begin{aligned}}
  {\end{aligned}\right\rbrace}


%

%





%

\pagestyle{headings} 

%


%

%

\newcommand{\Sym}{{\rm Sym}}

\newcommand{\vfi}{{\varphi}}
\newcommand{\Vfi}       {\bar \varphi}

 \def\Xint#1{\mathchoice
 {\XXint\displaystyle\textstyle{#1}}%
 {\XXint\textstyle\scriptstyle{#1}}%
 {\XXint\scriptstyle\scriptscriptstyle{#1}}%
 {\XXint\scriptscriptstyle\scriptscriptstyle{#1}}%
 \!\int}
 \def\XXint#1#2#3{{\setbox0=\hbox{$#1{#2#3}{\int}$}
 \vcenter{\hbox{$#2#3$}}\kern-.5\wd0}}
 \def\ddashint{\Xint=}
 \def\dashint{\Xint-}

\def\B{\mathbb{B}}
\def\test{\mathcal{D}}
\def\F{\mathbb{F}}
\def\N{\mathbb{N}}
\def\O{\Omega}
\def\R{\mathbb{R}}
\def\Kor{{\rm Kor}}
\def\pod#1{\mathop{#1}\limits}
\def\diagin{-\hskip-11.0truept\intop}
\def\diagint{{\raise-.1pt\hbox{--}\hskip-7.9pt\intop}}
\def\diagintop{\mathop{\mathchoice
{{\diagin}}%
{{\diagint}}%
{{\diagint}}%
{{\diagint}}%
}\limits}
 \def\ddashint{\Xint=}
 \def\dashint{\Xint-}

%
%
%
%

\section{Introduction}

This note provides existence and, more importantly, interior smoothness of solutions to a PDE system describing a static problem in a linearized yet nonlinear elasticity theory in an arbitrary space dimension. Let us begin with the model description. We consider an isotropic homogeneous elastic body occupying in a reference configuration a Lipschitz domain $\Omega \subset \er^d$. The body is affected by external forces of density $\bef:\Omega \to \er^d$ and surface forces $\bg:\Gamma_N \to \er^d$, where $\Gamma_N \subset \partial \Omega$. In addition, another $\Gamma_D$  part of the boundary $\partial \Omega$ is deformed by a displacement $\bu_0:\Gamma_D \to \er^d$. $\Gamma_N$ and $\Gamma_D$ are disjoint open subsets of $\partial \Omega$, whose union has the same $(d-1)$ Hausdorff measure as $\partial \Omega$. Then, the described deformation must satisfy the following balance of forces
\begin{equation}
\begin{aligned}
\label{eqa1}
- \Divv \bsig& = \bef  \quad (= -\Divv \bF) &&\textrm{in } \Omega, \\
\bsig \bn&= \bg \quad (=\bF \bn) &&\textrm{on }\Gamma_N,\\
\bu &= \bu_0 &&\textrm{on } \Gamma_D.
\end{aligned}
\end{equation}
Here, $\bsig:\Omega \to \er^{d\times d}_{sym}$ is the Cauchy stress tensor and $\bu:\Omega \to \er^d$ is the sought displacement field; $\bn$ is the outer normal unit vector. (For clarity, let us denote tensors by bold letters, vectors by regular letters and scalars by regular Greek letters. Accordingly, we disambiguate operators acting on tensors in bold and those acting on vectors in non-bold.) Already at this point, we assume that the external body and surface forces can be expressed by a given tensor $\bF:\Omega \to \er^{d\times d}_{sym}$. Naturally, if the problem is solvable, then any couple $(\bef, \bg)$ can be described (non-uniquely) by a tensor $\bF$. Additionally, in the case of the pure Neumann problem $\Gamma_N = \partial \Omega$, we avoid thereby using the necessary compatibility condition $\int_{\partial \Omega} g + \int_{\Omega} f = 0$, since it is already encoded into the existence of $\bF$.

In order to complete the problem \eqref{eqa1}, it remains to prescribe the constitutive relations for the Cauchy stress $\bsig$. The first classical law of linearized elasticity for isotropic homogeneous material is the generalized Hooke law, which has the form
\begin{equation}\label{Hook}
\bsig = 2\mu \, \D u + \lambda \divv\bu \;\bI = 2\mu\,  \Dd u + (2\mu/d + \lambda) \diver \bu \;\bI,
\end{equation}
where the constants $\mu$ and $\lambda$ are the so-called Lam\'{e} coefficients, the shear and the bulk modulus, respectively.  Here, we denoted by $\D u$ the linearized strain tensor and by $\Dd u$ its deviatoric part, i.e.,
$$
\D u:= \frac12 (\nablat \bu + (\nablat \bu)^T), \qquad \Dd u := \D u - \frac{\diver \bu}{d} \; \bI.
$$
Notice, that such a setting corresponds to a stored energy of the form
$$
W(\bu):= \mu |\Dd \bu|^2 + \tilde \lambda |\diver \bu|^2,
$$
where $ \tilde \lambda := \frac{2\mu+d\lambda}{2d}$. Consequently, natural constraints on the coefficients are $\mu>0$ and ${2\mu+d\lambda} >0$. Furthermore, the problem \eqref{eqa1} can be now equivalently restated as
\begin{equation}
\min_{\bu: \bu=\bu_0 \textrm{ on } \Gamma_D}\int_{\Omega} W(\bu) - \bF \cdot \D\bu \; dx.
\label{min}
\end{equation}
For suitably chosen data one can always find a unique weak solution $u$. In addition, due to linearity of the problem, it is then standard to show higher regularity properties of $u$ that are restricted only by smoothness of $\bF$, $\bu_0$ and $\partial \Omega$.

However, it is experimentally established that the Hooke law  \eqref{Hook} is not valid anymore for a body undergoing a large loading. Recently, many phenomenological laws allowing for more involved constitutive relations between the Cauchy stress and the linearized strain tensor were investigated. Although they are far from being linear, they can be theoretically justified even within the theory of small deformations (i.e., the `linearized' theory), see \cite{Ra07}. For an isotropic material, the most general constitutive law falling into the framework of the theory developed in \cite{Ra07} reads
\begin{equation}
\bsig = \alpha_1 \bI + \alpha_2 \D \bu + \alpha_3 \D\bu \, \D \bu, \label{Hook-b}
\end{equation}
where $\alpha_i$ may depend on all invariants of $\D \bu$. Decomposing further the dependence on the symmetric displacement gradient into its deviatoric and the trace part as well as requiring further that the material is hyperelastic, i.e., that there exists a potential - the generalized stored energy $W$, one arrives at the relation
\begin{equation}
\bsig = 2\mu(|\D^d \bu|) \D^d \bu + \tilde{\lambda}(\diver \bu) \diver \bu \bI, \label{Hook-c}
\end{equation}
where now $\mu$ and $\tilde{\lambda}$ are nonnegative scalar functions. The corresponding stored energy is then of the form
\begin{equation}
\label{stored-b}
W(\bu):= \psi(|\D^d \bu|) + \varphi(\diver \bu), \qquad \psi(s):= \int_{0}^s 2\mu(t)t\; dt, \qquad \varphi(s):=\int_0^s \tilde{\lambda}(t)t\; dt.
\end{equation}
Both functions $\psi$ and $\varphi$ are assumed to be nonnegative, convex and vanishing at zero. With this notation, solvability of \eqref{eqa1} supplemented with the constitutive relation \eqref{Hook-c} is still equivalent to minimization of \eqref{min} with the stored energy \eqref{stored-b}. Hence the solvability (in the weak sense) directly follows from the assumed convexity of $\psi$ and $\varphi$, provided reasonable growth and coercivity conditions for $W$ are assumed. Furthermore, an application of the standard difference quotient method enables to show certain regularity properties of the (unique) solution (at most the interior $W^{2,2}$ regularity for strongly coercive potentials). However, any more complex regularity theory, e.g. the $\mathcal{C}$-, $\mathcal{C}^{1}$- or $\mathcal{C}^{\infty}$-everywhere regularity is missing in general. The positive results in this direction for a rather general class of $W$'s in \eqref{stored-b} are known only either in the two dimensional setting or for general dimensional case if $\mu$ and $\lambda$ are `almost' constant functions, compare \cite{FrNe88} or \cite{KaMaSt99} for incompressible fluids.

In this short paper, we provide further regularity (smoothness) properties of a weak solution to \eqref{eqa1}, \eqref{Hook-c}  in a general multidimensional setting, that go much beyond the classical results. Indeed, we are able to cover certain highly nonlinear dependences of the Cauchy stress on the small strain tensor. It is worth noticing that the investigated problem falls into the class of `generalized' elliptic\footnote{We are using the word `generalized' here, because the potential $W$ depends only the symmetric gradient.} problems, where one cannot expect the full regularity of a solution, see the counterexamples in \cite{Ne77,SvYa02}. The only known structural assumption allowing for the full regularity reads
$$
W(\bu) \sim \tilde{W}(|\nabla \bu|),
$$
due to the classical result by Uhlenbeck \cite{Uh77}. However, for a nonlinear stored energy of the form \eqref{stored-b}, such a result is currently not available, due to `non-diagonality' of the corresponding elliptic operator and nonexistence of a proper substitute to the Uhlenbeck's method. Nevertheless, we show that for a certain class of models within \eqref{stored-b} one can overcome these difficulties and improve the smoothness of the solution significantly. Namely, our main result can be summarised as follows (for the fully rigorous formulation we refer the reader to Theorem~\ref{T:PR})
\begin{teo}\label{T:main}
Let $\bu_0, \bF$ be smooth and  $\mu$ be a constant. In \eqref{stored-b}, let $\psi(s):=\mu s^2$ and $\varphi$ be a smooth convex function, having at most the polynomial growth. Then there exists a unique weak solution $\bu$ to the problem \eqref{eqa1} belonging to $\mathcal{C}^{\infty}_{loc}(\Omega)$.
\end{teo}
Let us emphasize here, that this is the first regularity result for nonlinear systems od PDEs'  arising in the linearized elasticity theory, with no data-smallness or low-dimensionality restrictions. Despite being clearly far from covering the most natural case of \eqref{stored-b} (namely, both $\psi$ and  $\varphi$ being reasonably nonlinear), our result can be seen as the first step forward in the regularity theory for such problem. Beyond its theoretical novelty, it is viable applicatively: it covers real-world materials whose shear modulus remains constant, but the bulk modulus may change drastically with respect to the volume changes. Their examples are  rubber, certain polymers, concrete etc., see e.g. \cite{Hi91,Ho71,Ma92} and the references therein.

\section{Notation and definitions}

Firstly, let us introduce function spaces relevant to the stored energy \eqref{stored-b} with a constant $\mu$. Since we want to keep the allowed bulk modulus nonlinearity ($\varphi$) possibly general, we resort to Orlicz growths in our analysis. Nevertheless, let us accentuate at the very beginning that our results are new even for power--law growths. Let us briefly recall the Orlicz framework.  Since we intend to keep this part as concise as possible, we refer the interested reader e.g. to \cite[Appendix]{BurKap15} for more details. We shall start with the notion of the $\eN$ function.
\begin{mydef}[$\eN$-function]\label{def:app_orlicz:nf}
A real function $\vfi: \er \to \er_+$ is an \emph{$\eN$-function} iff it is even and there exists $\vfi'\!: \er_+ \to \er_+$
\begin{itemize}
\item[(N1)] \label{N1} that is right-continuous, non-decreasing,
\item[(N2)] \label{N2} that satisfies $\vfi' (0) =0$,  $ \vfi' (t) >0$ for $t>0$ and $\vfi' (+ \infty^-) =+ \infty$,
\end{itemize}
such that for all $t>0$ it holds
\begin{equation*}
\vfi (t) = \int_0^t \vfi' (s) \, ds.
\end{equation*}
\end{mydef}
Notice here, that it follows from (N\ref{N1}) that $\vfi$ is convex. Thus, we can also introduce its convex conjugate $\varphi^*$ by the formula $\varphi^*(s):=\sup_{t}(st-\varphi(t))$. Next,  let us define the Orlicz class $L^{\vfi}(\Omega)$ as
$$
L^{\vfi}(\Omega):=\left\{u\in L^1(\Omega): \; \int_{\Omega} \vfi(u) \dx <\infty \right\}.
$$
It becomes the Banach space for $\vfi$ satisfying the so--called $\Delta_2$ condition\footnote{For $\vfi$ not satisfying $\Delta_2$ condition, the relevant function space related to the corresponding Orlicz class is then defined as a union of such functions $u$ for which there exists $\lambda>0$ such that $\lambda u \in L^{\vfi}(\Omega)$.}, i.e.,
\begin{equation}\label{Delta2}
\vfi(2t)\le C(\vfi(t)+1), \quad \textrm{ for some } C>0 \textrm{ and all } t\in \mathbb{R}.
\end{equation}
The importance of the $\Delta_2$ condition appears also in the regularity theory, since it directly implies\footnote{Indeed, for $t\ge 0$ we can use the fact that $\varphi'$ is nondecreasing and nonnegative to observe that
$$
\varphi(t)=\int_0^t \varphi'(\tau) \le t\varphi'(t) \le \int_t^{2t}\varphi'(\tau)\le \varphi(2t) \overset{\eqref{Delta2}}{\le} C(\varphi(t)+1).
$$
For more details about the $\Delta_2$ condition and the above so--called good $\varphi'$ property see also \cite[Appendix]{BurKap15}, \cite{DieEtt08} or  \cite{DieKap13}.}
\begin{equation}
\varphi(t)\le t\varphi'(t) \le C(\varphi(t)+1).\label{delta2-better}
\end{equation}
Let us provide most typical growths that stay within our Orlicz structure and satisfy $\Delta_2$ condition.
Let $\kappa \ge 0$, $p > 1$. The classical power-law growths
\[
\vfi_1 (t) = \int_0^t \! (\kappa + s^{p-2}) \, s \,ds,  \qquad  \vfi_2 (t) = \int_0^t \!(\kappa + s^2)^\frac{p-2}{2} \, s \,ds
\]
are naturally allowed.

An example of an admissible $\eN$-function related to  a non-polynomial growth reads
\[
   \vfi_3 (t) = \int_0^t (\kappa +s^2)^\frac{p-2}{2}  s \ln^\beta (e+s) ds, \quad \beta >0.
   \]

 In what follows we will use standard notions of Lebesgue and Sobolev spaces $L^p$, $W^{k,p}$ respectively, where the subscript $_{\Gamma_D}$ will indicate that the considered functions vanish on $\Gamma_D$ (the Dirichlet part of the boundary). Next, let us introduce generalized Sobolev--Orlicz classes of vector--valued functions compatible with our problem setting. Here, $\bu_0:\Omega \to \mathbb{R}^d$ is a given measurable function.
$$
\begin{aligned}
WD_{\Gamma_D}&:=\left\{\bu\in W^{1,1}_{\Gamma_D}(\Omega; \mathbb{R}^d): \; \D^d\bu \in L^2(\Omega; \er^{d\times d}), \; \diver \bu \in L^{\vfi}(\Omega) \right\},\\
WD_{\bu_0}&:=\left\{\bu\in W^{1,1}_{\Gamma_D}(\Omega; \mathbb{R}^d): \; \bu=\bu_0 \textrm{ on } \Gamma_D, \; \D^d\bu \in L^2(\Omega; \er^{d\times d}), \; \diver \bu \in L^{\vfi}(\Omega) \right\}.
\end{aligned}
$$
Please notice that if $\vfi$ satisfies the $\Delta_2$ condition, $\Omega$ is Lipschitz  and $\bu_0$ has sufficiently highly integrable first order derivatives, then the notion $\bu \in WD_{u_0}$ is equivalent to $(\bu-\bu_0) \in WD_{\Gamma_D}$. Furthermore, in the case of $\Gamma_D = \emptyset$, we shall consider all functions belonging to $WD_{\Gamma_D}$ and $WD_{\bu_0}$ up to a rigid body motions, i.e., modulo all linear functions fulfilling $\D\bu \equiv 0$, in order to guarantee uniqueness of the desired displacement.

We shall frequently use a generic constant  $C>0$  that depends only on the data of our problem. It may generally vary line to line. If we need to trace any data dependences more precisely,  it will be clearly indicated in the text.

To conclude this section, let us introduce a notion of a weak solution. Below, we use the decomposition of (symmetric) $\bF$ into the deviatoric and the traceless part, i.e., $\bF=\bF^d + d^{-1} \tr \bF \, \bI$.
\begin{mydef}\label{D1}
Assume that $\Omega\subset \mathbb{R}^d$ is a Lipschitz domain and   $\bu_0 \in W^{1,2}(\Omega, \mathbb{R}^d)$ is such that $\diver \bu_0 \in L^{\varphi}$. Further, let $\bF$ be such that $\bF^d \in L^2(\Omega;\mathbb{R}^{d\times d}_{sym})$ and $\tr \bF \in L^{\varphi^*}(\Omega)$.  Let $\varphi$ be an $\eN$ function that satisfies \eqref{Delta2} and $\mu$ be a positive constant. We say that  $\bu \in  WD_{\bu_0}$ is a weak solution to \eqref{eqa1} with the constitutive law \eqref{Hook-c}--\eqref{stored-b} iff for all $\bv \in WD_{\Gamma_D}$ there holds
\begin{equation}\label{w-f}
\int_\Om 2\mu\, \D^d \bu \cdot \D\bv+  \varphi'(\diver \bu) \diver \bv  =  \int_\Om \bF \cdot  \D\bv.
\end{equation}
\end{mydef}
The Neumann part (formally) cancels out thanks to our `compatibility condition', i.e. the use of $\bF$ in \eqref{eqa1}. Observe that our weak formulation is meaningful. Indeed, the critical term can be estimated thanks to (N\ref{N1}) and \eqref{delta2-better} as follows
$$
|\varphi'(\diver \bu) \diver \bv| \le |\varphi'(\diver \bv) \diver \bv| + |\varphi'(\diver \bu) \diver \bu| \le C(1+ \varphi(\diver \bu) + \varphi(\diver \bv)) \in L^1(\Omega).
$$

\section{Result}\label{sec:rel}
\begin{teo}\label{T:PR}
Let all assumptions of Definition~\ref{D1} be satisfied. Then, there exists a weak solution to \eqref{eqa1} such that
\begin{align}
\|\D^d\bu\|^2_{L^2(\Omega)} + \int_{\Omega}\varphi(\diver \bu)&\le C(\mu,d)\left(\|\bF^d\|^2_{L^2(\Omega)}+\|\D^d\bu_0\|^2_{L^2(\Omega)} + \int_{\Omega}\varphi^*(\tr \bF)+\varphi(\diver \bu_0)\right)=:A, \label{TE1}
\end{align}
This weak solution is unique, provided that $\varphi$ is strictly convex. More surprisingly, it enjoys the following smoothness properties in any compact set $K\subset \Omega$:

\begin{align}
\left\|\frac{\partial \bu_i}{\partial x_j} - \frac{\partial \bu_j}{\partial x_i}\right\|_{W^{k,p}(K)}+\left\|\frac{2\mu(d-1)}{d}\diver \bu  + \varphi'(\diver \bu)\right\|_{W^{k,p}(K)} &\le C(K,\Omega, A,k,p)(\|\bF\|_{W^{k,p}(\Omega)}+1), \label{TE2}\\
\|\nabla \bu\|_{W^{k,p}(K)}&\le C(K, \Omega, A,k,p, \|\varphi\|_{\mathcal{C}^{k+2}}, \|\bF\|_{W^{k,p}(\Omega)}),\label{TE3}
\end{align}
where \eqref{TE2} is valid for all $k\in \mathbb{N}_0$ and all $p\in (1,\infty)$, whereas \eqref{TE3} holds both for $k=0,1$ with all $p\in (1,\infty)$  and for arbitrary $k\in \mathbb{N}$ with any $p\in (d/2, \infty)$.
\end{teo}

\section{Proof of Theorem~\ref{T:PR}}\label{sec:proof}

\paragraph{Weak existence \& uniqueness}
We use the direct methods of  the calculus of variations to obtain the existence. Indeed, let us mimic the problem \eqref{min} and can seek for a minimizer $\bu\in WD_{\bu_0}$ fulfilling for all $\bv \in WD_{\bu_0}$
\begin{equation} \label{min-g}
\int_{\Omega}\mu|\D^d \bu|^2 + \varphi(\diver \bu)  - \bF \cdot \D \bu \le \int_{\Omega}\mu|\D^d \bv|^2 + \varphi(\diver \bv)  - \bF \cdot \D \bv.
\end{equation}
The functional is definitely coercive and due to convexity of $\varphi$ and its superlinearity at infinity, we can find a minimizer fulfilling \eqref{min-g}. In addition, the uniform estimate \eqref{TE1} easily follows. Furthermore, since $\varphi$ satisfies the $\Delta_2$ condition, and consequently $\bv:=\bu+\bw$ with an arbitrary $\bw \in WD_{\Gamma_0}$ is an admissible competitor in \eqref{min-g},  we can easily derive the Euler--Lagrange equations for \eqref{min-g}, which is nothing else than the identity \eqref{w-f}. The uniqueness then follows from the strict convexity of $\varphi$. \qed

We are now approaching the crucial part, namely our proof of the interior smoothness. Before passing to details, let us observe that the main idea behind this proof is the following observation: Since
\[
\diver \Divv \Dd f = \frac{d-1}{d} \Delta \diver f,
\]
then $\diver \eqref{eqa1}$ produces a scalar, well-manageable equation for $\diver u$.

\paragraph{Estimate \eqref{TE2}}
We start with deriving an elliptic equation for $\diver \bu$.
Let $v\in \mathcal{C}^2_0(\Omega)$ be arbitrary. Setting $\bv:=\nabla v$ in \eqref{w-f} we observe that
\begin{equation}\label{MB1}
\int_\Om 2\mu \D^d \bu \cdot \nabla^2 v +  \varphi'(\diver \bu) \Delta v  =  \int_\Om \bF \cdot  \nabla^2 v.
\end{equation}
Next, let $G$ be the Green function to the Laplace equation in  $\er^d$ and let us define (in the sense of distribution)
$$
g:=G*\diver \Divv \bF,
$$
where we extend  $\bF\equiv 0$ outside $\Omega$.
Notice that such a $g$ solves the problem
\begin{equation}\label{dfg}
\int_{\Omega} g \Delta v = \int_{\Omega} \bF \cdot \nabla^2 v \qquad \textrm{ for all } v \in \mathcal{C}^{\infty}_{0}(\Omega).
\end{equation}
In addition for any compact $K \subset \Omega$, all $k\in \mathbb{N}$ and all $p\in (1,\infty)$, we have the estimate
\begin{equation}
\|g\|_{W^{k,p}(K)} \le C(K,k,p)\|\bF\|_{W^{k,p}(\Omega)}. \label{g-est}
\end{equation}

Finally, using integration by parts, it is not difficult to deduce that
$$
\begin{aligned}
\int_\Om \D^d \bu \cdot \nabla^2 v=\int_\Om \D^d \bu \cdot (\nabla^2 v - d^{-1}\Delta v \, \bI)=\int_\Om \nabla \bu \cdot (\nabla^2 v - d^{-1}\Delta v \, \bI)=\frac{d-1}{d}\int_\Om \diver \bu \Delta v.
\end{aligned}
$$
Consequently, using also \eqref{MB1} and \eqref{dfg}, we obtain
$$
\int_{\Omega} \left(\frac{2\mu(d-1)}{d}\diver \bu  + \varphi'(\diver \bu) -g\right)\Delta v =0,
$$
which means that $\frac{2\mu(d-1)}{d}\diver \bu  + \varphi'(\diver \bu) -g$ is harmonic in $\Omega$. Therefore, for any compact $K\subset \Omega$ we have that
$$
\left\|\frac{2\mu(d-1)}{d}\diver \bu  + \varphi'(\diver \bu) -g\right\|_{W^{k,p}(K)} \le C(K)\left\|\frac{2\mu(d-1)}{d}\diver \bu  + \varphi'(\diver \bu) -g\right\|_{L^1(\Omega)}\le C(A,K),
$$
where $A$ is defined in \eqref{TE1}. Thus, using also \eqref{g-est}, we finally deduce that
\begin{equation}
\left\|\frac{2\mu(d-1)}{d}\diver \bu  + \varphi'(\diver \bu)\right\|_{W^{k,p}(K)} \le C(K,A,p,k)(\|\bF\|_{W^{k,p}(\Omega)}+ \|\diver \bu\|_{L^\varphi (\Omega)}) \label{fin-div}.
\end{equation}
Hence the second part of the estimate \eqref{TE2} follows from estimating the right hand side of \eqref{fin-div} with the help of \eqref{TE1}.

Similarly, using \eqref{w-f}, we deduce that in the distributional sense
\begin{equation}\label{w-f-c}
-\mu \Delta \bu + \mu \nabla \diver \bu - \nabla \left( \frac{2\mu(d-1)}{d}\diver \bu  + \varphi'(\diver \bu)\right) = -\Divv \bF.
\end{equation}
Consequently, taking  distributionally aa $i,j$ component of $\nabla \wedge$ in \eqref{w-f-c} yields for any $i,j=1,\ldots,d$
\begin{equation}\label{curl}
\mu \Delta \left (\frac{\partial \bu_i}{\partial x_j} - \frac{\partial \bu_j}{\partial x_i}\right ) = \sum_{k=1}^d \frac{\partial}{\partial x_k} \left(\frac{\partial \bF_{ik}}{\partial x_j} - \frac{\partial \bF_{jk}}{\partial x_i}\right).
\end{equation}
Thus, since we already know that $\nabla \bu \in L^{1}$, we can use the singular integral theory to conclude that for all $i,j=1,\ldots,d$, all $k\in \mathbb{N}$, all $p\in (1,\infty)$ and all compact set $K\subset \Omega$
\begin{equation}\label{est}
\left\|\frac{\partial \bu_i}{\partial x_j} - \frac{\partial \bu_j}{\partial x_i}\right\|_{W^{k,p}(K)}\le C(\mu, k,p,K)(\| \nabla \bu\|_{L^{1}(\Omega)}+\|\bF\|_{W^{k,p}(\Omega)})
\end{equation}
and \eqref{TE2} follows from \eqref{est} and \eqref{TE1}.\qed

\paragraph{Estimate \eqref{TE3}}
Going back to \eqref{w-f-c} and using the theory for the Laplace equation, we directly deduce the estimate
\begin{equation}\label{26}
\|\bu\|_{W^{k+1,p}(K')} \le C(k,p,K,K')\left(\|\bu\|_{L^2(\Omega)} +\|\bF\|_{W^{k,p}(\Omega)} + \|\diver \bu\|_{W^{k,p}(K)} + \left\|\frac{2\mu(d-1)}{d}\diver \bu  + \varphi'(\diver \bu)\right\|_{W^{k,p}(K)} \right)
\end{equation}
for arbitrary compact sets $K'$, $K$ such that there is an open set $\Omega'$ fulfilling $K'\subset \Omega'\subset K\subset \Omega$, arbitrary $k\in \mathbb{N}_0$ and arbitrary $p\in (1,\infty)$. Therefore due to the already obtained estimate \eqref{fin-div}, we see that it is enough to get the bound on $\diver \bu$ which will be however read again from \eqref{fin-div}.

We start with the case $k=0$. Since $\varphi'(s)s\ge 0$ (which follows from the fact that $\varphi$ is even), hence $|s| \le |\varphi'(s)+s|$ and one has
$$
\|\diver \bu\|_p \le C(d,\mu)\left\|\frac{2\mu(d-1)}{d}\diver \bu  + \varphi'(\diver \bu)\right\|_{p}.
$$
Consequently, it follows from \eqref{26} and \eqref{fin-div} that
\begin{equation}\label{27}
\|\bu\|_{W^{1,p}(K')} \le C(k,p,K,K')\left(\|\bu\|_{L^2(\Omega)} +\|\bF\|_{L^{p}(\Omega)} \right)
\end{equation}
and \eqref{TE3} holds for $k=0$ and $p\in (1,\infty)$.
Secondly, take $k=1$. Since $\varphi'$ is nondecreasing (and Lipschitz, recall  \eqref{TE3} for $k=1$), it is easy to observe that
$$
\|\diver \bu\|_{1,p} \le C(d,\mu)\left\|\frac{2\mu(d-1)}{d}\diver \bu  + \varphi'(\diver \bu)\right\|_{1,p}
$$
for all $p\in (1,\infty)$. Consequently, we again get that (now we set $k=1$ in \eqref{26})
\begin{equation}\label{28}
\|\bu\|_{W^{2,p}(K')} \le C(k,p,K,K')\left(\|\bu\|_{L^2(\Omega)} +\|\bF\|_{W^{1,p}(\Omega)} \right)
\end{equation}
and \eqref{TE3} for $k=1$ again follows.
Finally, for the remaining range of $k$'s we shall restrict ourselves to the case $p> d/2$. From \eqref{fin-div} one obtains
\begin{equation}
\left\|\left(\frac{2\mu(d-1)}{d}+\varphi''(\diver \bu)\right) \nabla \diver \bu\right\|_{W^{k,p}(K')} \le C(K')(\|\bF\|_{W^{k+1,p}(\Omega)}+ \|\diver \bu\|_{L^\varphi(\Omega)}). \label{fin-div2}
\end{equation}
From \eqref{28} and using $p> d/2$
\begin{equation}\label{eq:i1}
\|\diver \bu\|_{L^\infty(K')} \le C(d,p,K',\Omega)\left(\|\bu\|_{L^2(\Omega)} +\|\bF\|_{W^{1,d+\delta}(\Omega)} \right) \le C(d,p,K',\Omega)\left(\|\bu\|_{L^2(\Omega)} +\|\bF\|_{W^{k+1,p}(\Omega)} \right),
\end{equation}
therefore arguments of $\varphi''$ and its derivatives will remain bounded.
Let us now show inductively that for $p> d/2$ there holds
\begin{equation}
\left\|\nabla \diver \bu\right\|_{W^{k,p}(K')} \le C(K', \Omega, d, p,\|\varphi\|_{\mathcal{C}^{k+2}}, k, A) (1+ \|\bF\|_{W^{k+1,p}(\Omega)})^{3^k}. \label{fin-div3}
\end{equation}
This relation is valid for $k=0$ thanks to \eqref{fin-div2} and  \eqref{TE1}. Assuming now that \eqref{fin-div3} holds for $k-1$, we shall show that it also holds for $k$. One has

\begin{equation}\label{eq:ind}
\left\|\nabla \diver \bu\right\|_{W^{k,p}(K')} =\left\|\nabla^{k}  \nabla \diver \bu\right\|_{L^{p}(K')} + \left\|\nabla \diver \bu\right\|_{W^{k-1,p}(K')} := I + II.
\end{equation}
Furthermore
\begin{equation}\label{eq:i2}
\begin{aligned}
I &\le \frac{d}{2\mu(d-1)}  \left\| \left(\frac{2\mu(d-1)}{d}+\varphi''(\diver \bu)\right)  \nabla^k \left( \nabla \diver \bu \right)\right\|_{L^{p}(K')}  \\
& \le C(\mu, d) \left\| \nabla^k \left(\left(\frac{2\mu(d-1)}{d}+\varphi''(\diver \bu)\right) \nabla \diver \bu \right)\right\|_{L^{p}(K')} +  \sum^k_{i=1} C  \left\|  \nabla^i ( \varphi''(\diver \bu) )\nabla^{k-i}\nabla  \diver \bu \right\|_{L^{p}(K')}  \\
& \le C(\mu, d) \left\| \nabla^k \left(\left(\frac{2\mu(d-1)}{d}+\varphi''(\diver \bu)\right) \nabla \diver \bu \right)\right\|_{L^{p}(K')} + C  \|\varphi\|_{\mathcal{C}^{k+2}}  \|\diver \bu\|_{L^\infty(K')} \|\diver \bu\|^{{2}}_{W^{k,2p}(K')}.
\end{aligned}
\end{equation}
The last term of \eqref{eq:i2}, via the inductive assumption \eqref{fin-div3} for $k-1$,  \eqref{eq:i1} and next by $W^{k+1,p}\hookrightarrow W^{k,2p}$ for $p> d/2$ can be estimated with
\[
 C \left(\|\bu\|_{L^2(\Omega)} +\|\bF\|_{W^{k+1,p}(\Omega)} \right) (\|\bF\|_{W^{k,2p}(\Omega)}+ \|\diver \bu\|_{L^\varphi(\Omega)})^{2 \cdot 3^{k-1}} \le  C  ( 1+ \|\bF\|_{W^{k+1,p}(\Omega)})^{3^k},
\]
where we have used \eqref{TE1} and $C= C(K', \Omega, d, p,\|\varphi\|_{\mathcal{C}^{k+2}}, k, A)$.
Since the last but one  of \eqref{eq:i2} is controlled by \eqref{fin-div2}, it holds
\[
I \le C ( 1+ \|\bF\|_{W^{k+1,p}(\Omega)})^{3^k},
\]
Since  the lower-order term $II$ is estimated by the inductive assumption \eqref{fin-div3} for $k-1$, we have arrived at validity of  \eqref{fin-div3} for $k$.
Estimate \eqref{TE3} is hence proven. \qed

\bibliographystyle{abbrv}
\bibliography{references}
\end{document}